\documentstyle[amsfonts,amssymb,12pt]{article}
\input xypic
\input cyracc.def

\newfam\cyifam
\font\tencyi=wncyi10
\font\sevencyi=wncyi7
\font\fivecyi=wncyi5

\textfont\cyifam=\tencyi \scriptfont\cyifam=\sevencyi
  \scriptscriptfont\cyifam=\fivecyi
\newfam\cyrfam
\font\tencyr=wncyr10
\font\sevencyr=wncyr7
\font\fivecyr=wncyr5
\def\cyr{\fam\cyrfam\tencyr\cyracc}
\textfont\cyrfam=\tencyr \scriptfont\cyrfam=\sevencyr
  \scriptscriptfont\cyrfam=\fivecyr

\renewcommand{\paragraph}{}

\newcommand{\ad}{{\mathrm ad}}
\newcommand{\rar}{\rightarrow}

\newcommand{\Proof}{{\bf Proof}.\, }

\newcommand{\fe}{{\cal O}^{\mathrm pre}}

\newcommand{\Def}{\mathsf{Def}}
\newcommand{\Df}{\mbox{\cyr Def}}
\newcommand{\lon}{\longrightarrow}

\newcommand{\Z}{{\Bbb Z}}
\newcommand{\p}{{\partial}}

\newcommand{\C}{{\Bbb C}}
\newcommand{\R}{{\Bbb R}}
\newcommand{\ot}{\otimes}
\newcommand{\tl}{\tilde}

\newcommand{\Beq}{\begin{equation}}
\newcommand{\Eeq}{\end{equation}}
\newcommand{\Beqr}{\begin{eqnarray}}
\newcommand{\Eeqr}{\end{eqnarray}}
\newcommand{\Beqrn}{\begin{eqnarray*}}
\newcommand{\Eeqrn}{\end{eqnarray*}}
\newcommand{\Ba}{\begin{array}}
\newcommand{\Ea}{\end{array}}
\newcommand{\Bi}{\begin{itemize}}
\newcommand{\Ei}{\end{itemize}}
\newcommand{\Bc}{\begin{center}}
\newcommand{\Ec}{\end{center}}
\newcommand{\fg}{{\frak g}}

\newcommand{\f}{{\cal O}}

\newcommand{\cM}{{\cal M}}
\newcommand{\cA}{{\cal A}}

\newcommand{\cU}{{\cal U}}

\newcommand{\cV}{{\cal V}}
\renewcommand{\cH}{{\cal H}}
\newcommand{\cX}{{\cal X}}
\newcommand{\cS}{{\cal S}}

\newcommand{\al}{\alpha}
\newcommand{\be}{\beta}

\newcommand{\Ga}{\Gamma}

\newcommand{\var}{\varepsilon}

\newcommand{\dal}{\dot{\alpha}}
\newcommand{\dbe}{\dot{\beta}}

\newcommand{\tlo}{\tilde{0}}
\newcommand{\tln}{\tilde{1}}

\newcommand{\bp}{\bar{\partial}}

\newcommand{\Ker}{{\mathsf Ker}\, }
\newcommand{\Img}{{\mathsf Im}\, }

\newcommand{\Hom}{\mbox{Hom}}
\newcommand{\bH}{{\bf H}}

\newcommand{\cT}{{\cal T}}

\newcommand{\sip}{\smallskip}
\newcommand{\bip}{\bigskip}

\begin{document}

\title{A note on extended complex manifolds}
\author{ S.A.\ Merkulov}
\date{}
\maketitle

\sloppy

\begin{abstract}
We introduce a category of extended complex manifolds, and prove
that the functor describing deformations of a classical compact
complex manifold $M$ within this category is versally representable by
(an analytic subspace  in) $\bH^*(M,T_M)$.

By restricting the associated
versal family of extended complex manifolds over $\bH^*(M,T_M)$
to the subspace $\bH^1(M,T_M)$ one gets a correct limit
to the  classical picture.

\end{abstract}

\bip

\bip

{\bf 0. Introduction.} Recent advances in homological mirror symmetry
\cite{Ba,BaKo,Ko1,Ma1,Ma2,P,PZ} give a strong evidence for the existence
 of {\em extended}\,  versions of at least three classical categories in 
algebraic geometry ---
the categories of compact 
complex manifolds, of coherent sheaves, and of submanifolds. 
One such category
(of  extended special Lagrangian submanifolods)
has been constructed in \cite{Me1}. In this note we suggest an extension
of the category of complex manifolds.

\sip

Curiously, even with no appropriate notion of extended complex structure
at hand, one 
can nevertheless say a lot about their moduli spaces  \cite{Ba,BaKo,Ma1,Me2}.
This egg-before-chicken situation is due to a remarkably effective,
purely algebraic paradigm of the modern deformation theory (see \cite{GM,Ko2}
and references therein):
\Bi
\item Given a mathematical structure one wishes to deform,  the first step
should be a search for a differential $\Z$-graded (dg-, for short)
Lie $k$-algebra $(\fg=\bigoplus_{i\in \Z} \fg^i, d, [\ ,\  ])$ which
``controls'' the deformations. Next one defines a deformation functor
$$
\Ba{rccc}
{\Def}^{\, 0}_{\fg}: & \left\{\Ba{l} \mbox{the category of Artin}\\
                    \mbox{$k$-local algebras} \Ea \right\}&
\lon & \left\{\mbox{the category of sets}\right\}
\Ea
$$
as follows
$$
{\Def}^{0}_{\fg}(A)=\left\{\Gamma\in (\fg\ot m_{A})^1 \mid d\Gamma
+ \frac{1}{2}[\Gamma, \Gamma]=0\right\}/\exp{(\fg\ot m_{A})^0},
$$
where  $m_{A}$ is the maximal ideal of the Artin algebra  $A$,
the latter is viewed as a $\Z$-graded algebra concentrated in degree 
zero (so that $(\fg\ot m_{\cA})^i=\fg^i\ot m_{A}$), and the quotient 
is taken with respect to the following representation of the gauge group 
$\exp{(\fg\ot m_{A})^0}$, 
$$
\Gamma \rar \Gamma^g = e^{\ad_g}\Gamma - \frac{e^{ad_g}-1}{\ad_g}dg,
 \ \ \ \ g\in (\fg\ot m_{A})^0,
$$
where $\ad$ is just the usual internal automorphism of $\fg$,
$\ad_g \Gamma:= [g, \Gamma]$. Finally one tries 
to represent the deformation 
functor by a topological (pro-Artin) algebra $\f_{S}$ so that
$$
{\Def}^{\, 0}_{\fg}(A)= \Hom_{\mathrm cont}(\f_{S}, A).
$$
This associates to the given mathematical structure  a formal 
 moduli space $S$
whose ``ring of functions'' is $\f_{S}$.
\Ei

\sip

The tangent space, ${\Def}_{\fg}^{\, 0}(k[\var]/\var^2)$, to the functor
${\Def}_{\fg}^{\, 0}$ is isomorphic to the  cohomology group 
$\bH^1(\fg)$ of the complex $(\fg,d)$. If one extends in the obvious 
way the above deformation functor to the category of arbitrary 
$\Z$-graded $k$-local Artin algebras (which may not be concentrated 
in degree 0), one gets the functor ${\Def}_{\fg}^{\, *}$
with the tangent space isomorphic to  the full cohomology group
$\bH^*(\fg)=\oplus _{i\in \Z}\bH^i(\fg)$.

\sip

The dg-Lie algebra controlling  deformations of a given complex structure 
on a compact manifold $M$ is $(\fg:= \bigoplus_{i}
\Gamma(M, T_M\ot \Omega^{0,i}_M), \bp)$, where $T_M$ is the tangent holomorphic
sheaf, and $\Omega^{0,i}_M$ the sheaf of $(0,i)$ forms. 
The classical deformation functor
${\Def}_{\fg}^{\, 0}$ is known to be versally representable by the Kuranishi
analytic subspace in $\bH^1(\fg)=\bH^1(M,T_M)$. Its extension, ${\Def}_{\fg}^{\, *}$,
is the main technical tool for introducing and studying the moduli
space of so called {\em 
extended complex structures}\, \cite{BaKo}. Actually, the authors
of \cite{BaKo} go even further and study the 
extended deformation functor
associated with the larger dg-Lie algebra,  
$(\hat{\fg}=\Gamma(M, \wedge^{*}T_M\ot \Omega^{0,*}_M), \bp)$. 
They prove that ${\Def}_{\hat{\fg}}^{\, *}$ is non-obstructed provided $M$ is
a Calabi-Yau manifold,
and show that the associated extended moduli space, 
$S_{\mathrm ext}\simeq \bH^*(M,\wedge^{*}T_M)$, has an induced structure of 
Frobenius manifold \cite{BaKo,Ba}. For arbitrary $M$, the moduli
space $S_{\mathrm ext}$ is canonically an $F_{\infty}$-manifold \cite{Me2}.

\sip

There is, however, an obvious
problem with the above approach as it offers no
geometric explanation of what this
{\em extended complex structure}\, might be. 
It is an urgent and important 
problem \cite{Ma2} to find a 
geometric embodiment of this notion, develop its deformation theory,
and check that the base of the resulting versal deformation
can be canonically identified with what one gets from the above purely
algebraic approach. The present note offeres {\em a}\, realization of 
this programme for the extended deformation functor ${\Def}_{\fg}^{\, *}$.
It is unlikely, however, that the same approach will yield a 
geometrical model for 
the functor ${\Def}_{\hat{\fg}}^{\, *}$ --- this
 seems to be a much more intricate object.

\bip

{\bf 1. Extended complex manifolds.} 
These will be defined in
two steps. First comes the notion of pre-complex manifold.

\bip

{\bf 1.1. Model  pre-complex structures.} Let $U$ be an open domain in
$\R^{2n|0}$. 
A {\em pre-complex}\, structure on $U$ is a map
$$
\phi: U\lon \C^{2n|n}
$$
such that
\Bi
\item[(i)] $\phi$ is a smooth embedding, and
\item[(ii)] at each $x\in U$, the completions
of the associated stalks of the sheaves
$C^{\infty}(U)$ and $\phi^{-1}(\f_{\C^{2n|n}}/{\mathsf nilpotents})$
are isomorphic.
\Ei
Here $C^{\infty}(U)$ stands for the sheaf of  smooth 
complex valued functions on $U$, and $\f_{\C^{2n|n}}$  for the  sheaf of 
 holomorphic functions on the supermanifold ${\C^{2n|n}}$.
The ringed space 
$$
\cU^n =\left(U, \fe_U:= \phi^{-1}(\f_{\C^{2n|n}})\right)
$$
is called a {\em pre-complex domain of dimension}\, $n$.
The reduced sheaf $\fe_U/{\mathsf nilpotents}$ is  denoted sometimes by
$\fe_{U,{\mathrm rd}}$.

\bip

A {\em morphism}\, of pre-complex domains $\cU^n \rar \cV^m$ is, 
by definition, a pair $(f_{C^{\infty}}, [f])$, consisting of a smooth map 
$f_{C^{\infty}}: U\rar V$
and a  germ, $[f]$, of holomorphic maps $f$  from a small open neighbourhood
of $\Img(\phi_U)$ in $\C^{2n|n}$ to a small open neighbourhood 
of $\Img(\phi_V)$
in $\C^{2m|m}$ such that, for each $x\in U$,
$$
f_{C^{\infty}}^*\left(\overline{\overline{\fe_{V,{\mathrm rd},x}}}\right) = 
\overline{\overline{f^*(\fe_{V,{\mathrm rd},x})}}, 
$$
where\, $\overline{\overline{\hspace{3mm}^{}}}$\,
 stands for the  completion of the stalks at $x$
with respect to the natural ideals. We often abbreviate $(f_{C^{\infty}}, 
[f])$ simply to $f$.

\bip
{\bf 1.2. Fact.} Any open domain $U\subset \R^{2n|0}$ admits a pre-complex
structure. Indeed, identifying $\R^{2n|0}$ with $\C^{n|0}$ one 
immediately gets a required map $\phi$ by further identifying
 $C^{n|0}$ with the ``diagonal'' subspace in $\C^{2n|n}$ given, in the 
natural coordinates
$(z^{\al}, z^{\dal}, \psi^{\dal})$, $\al, \dal=1, \ldots,n$, on $\C^{2n|n}$,
by the equations
$$
z^{\dal}= \overline{(z^{\al})},\ \ \  \psi^{\dal}=0,
$$
where the bar denotes complex conjugation. 

\bip

{\bf 1.3. Definition.} An {\em n-dimensional  pre-complex
manifold}\  is a ringed space $\cM=(M,\fe_M)$ 
modeled on pre-complex domains of dimension $n$. 
It is called {\em compact}\, if the underlying
smooth manifold $M$ is compact.

\sip

A {\em holomorphic}\, 
vector field $v$ on $\cM$ is, as usual, a $\C$-linear
automorphism of the structure sheaf,
$$
v: \fe_M \lon \fe_M,
$$
whose restriction to each stalk, $\fe_{M,x}$, is a derivation
of the ring $\fe_{M,x}$. The sheaf of holomorphic vector fields
on $\cM$ is denoted by $\cT_{\cM}$.

\bip

{\bf 1.4. Definition.} An {\em extended complex manifold}\, is a pair
$(\cM, \p_{\cM})$ consisting of a   compact pre-complex manifold $\cM$
and an odd holomorphic vector field $\p_{\cM}$ such that
$[\p_{\cM}, \p_{\cM}]=0$\footnote{Odd vector fields with this 
property are often called {\em homological}.}
and $\dim 
\bH(\cM,\p_{\cM}) <\infty$. Here
$$
\bH(\cM,\p_{\cM}):= \frac{\Gamma(M, \Ker [\p_{\cM}, \ldots ])}{\Gamma(M,
\Img [\p_{\cM},\ldots])},
$$
with the differential  given by
$$
\Ba{rccc}
[\p_{\cM},\ldots]: &  \cT_{\cM} & \lon & \cT_{\cM} \\
& V & \lon & [\p_{\cM}, V].
\Ea
$$
Note that $\bH(\cM,\p_{\cM})$ is canonically a Lie superalgebra with the 
brackets, $[\ ,\ ]_{\bH}$, induced from the usual commutator of holomorphic 
vector fields.

\sip

One can also associate with $(\cM, \p_{\cM})$
 the {\em  cohomology manifold}, ${\cH}(\cM, \p_{\cM})$,
which is, by definition, a ringed space $(M, {\mathsf Ker}\, \p_{\cM}/
\Img \p_{\cM})$.

\sip

A {\em morphism}\, of extended complex manifolds,
$$
f: (\cM, \p_{\cM}) \lon (\cS, \p_{\cS}),
$$ 
is a morphism of the associated pre-complex manifolds, 
$f: (M,\fe_M)\rar (S,\fe_M)$,
 which commutes with the homological vector fields, i.e.
$$
\p_{\cM}\left(f^*(g)\right) = f^*(\p_{\cS}g)
$$
for any $g\in \fe_S$. One may reformulate this as $f_*(\p_{\cM})=\p_{\cS}$.

\bip

{\bf 1.5. Basic example.} Let $M$ be a compact complex manifold.
There is associated an extended complex manifold $(\cM=(M, \fe_{M}), 
\bar{\p})$ constructed as follows:
\Bi
\item[(i)] First consider a natural embedding,
$$
\Ba{rccc}
\phi: & M & \lon & M \times \Pi {T}_{\overline{M}} \\
&        x  & \lon & ({x}, O(x))

\Ea
$$
where ${T}_{\overline{M}}$ is the total space of the holomorphic
vector bundle over the conjugate complex manifold $\overline{M}$,
$\Pi$  the parity change functor, and 
$O:\overline{M} \rar \Pi {T}_{\overline{M}}$  the zero section.
\item[(ii)] Set  
$$
\fe_{M}:=\phi^{-1}\left({\mathsf the\ structure\ sheaf\ on}\
 M \times \Pi {T}_{\overline{M}}\right).
$$ 
\item[(iii)] Note that $\Pi {T}_{\overline{M}}$ is canonically
a dg-manifold (of dimension $n|n$) with the  homological field
being just the $(1,0)$-part of the de Rham differential on $\overline{M}$.
The latter induces a homological vector field on $\cM$ which we denote by 
$\bar{\p}$.
\Ei

\bip

{\bf 1.5.1. Remark.} Let $T_{M}$ be the holomorphic tangent sheaf
on $M$. 
 We make the cohomology $\bH^*(M, T_{M})=\bigoplus_k 
\bH^k(M, T_{M})$ into a superspace
by setting
$$
\bH^*(M, T_{M})_{\tlo}= \bigoplus_{k\ {\mathsf is\ odd}}
\bH^k(M, T_{M}),
\ \ \ \ \ 
\bH^*(M, T_{M})_{\tln}= \bigoplus_{k\ {\mathsf is\ even}}
\bH^k(M, T_{M}).
$$
This choice of $\Z_2$-grading is in agreement with the classical 
deformation theory
where $\bH^1(M, T_{M})$ is  even. Then the natural map,
$$
\Ba{ccc}
T_{M}\ot \Omega^{0,p}_{M}\ \times\ T_{M}\ot \Omega^{0,q}_{M}
&\lon & T_{M}\ot \Omega^{0,p+q}_{M} \vspace{2mm}\\
X\ot\al\ \ \  \times\ \ \ Y\ot \be & \lon & [X,Y]\ot(\al\wedge \be)
\Ea
$$
induces on $\bH^*(M, T_{M})$ the structure of an {\em odd}\, Lie
superalgebra (cf.\ \cite{Ma1}). Reversing the parity, we make
 $\Pi \bH^*(M, T_{M})$ into a Lie superalgebra.

\bip

{\bf 1.5.2. Lemma.} {\em Let $M$ be a complex manifold, and
$(\cM, \bar{\p})$ the associated extended complex manifold. Then
\Bi
\item[(a)] $\cH(\cM, \bar{\p})$ is precisely $M$ with its original 
complex structure.
\item[(b)] $\bH(\cM, \bar{\p})= \Pi \bH^*(M, T_{M})$ as 
Lie superalgebras.
\Ei
}

\sip

\Proof The statement (a) follows immediately from the Poincare
$\bar{\p}$-lemma. The second statement requires a small computation.
Let $\{z^{\al}\}$ be a  local coordinate system on $M$, and
$\{z^{\dal}, \psi^{\dal}=dz^{\dal}\}$ the associated local coordinate
system on $\Pi T_{\overline{M}}$. The collection $\{z^{\al}, 
z^{\dal}, \psi^{\dal}\}$ gives rise to a coordinate chart on $(M,\fe_{M})$
so that any $V\in \Gamma(M, \cT_{\cM})$ can be locally represented as
$$
V= \sum_{\al=1}^n V^{\al}\frac{\p}{\p z^{\al}} + \sum_{\dal=1}^n\left(V^{\dal}
\frac{\p}{\p z^{\dal}} + W^{\dal}\frac{\p}{\p \psi^{\dal}}\right),
$$
for some local sections, $V^{\al}, V^{\dal}, W^{\dal}$, of $\fe_{M}$.
As 
$$
\bar{\p}= \sum_{\dal=1}^n \psi^{\dal} \frac{\p}{\p z^{\dal}}, 
$$
we have
$$
[\bar{\p}_0, V]= \sum_{\al=1}^n (\bar{\p} V^{\al})
\frac{\p}{\p z^{\al}} + \sum_{\dal=1}^n \left(\bar{\p} V^{\dal} - 
(-1)^{\tilde{V}}W^{\dal}\right)\frac{\p}{\p z^{\dal}} +
\sum_{\dal=1}^n (\bar{\p} W^{\dal})\frac{\p}{\p \psi^{\dal}}.
$$
Thus $V\in \Gamma(M, \Ker [\bar{\p}, \ldots])$ if and only if
$\bp V^{\al}=0$, $W^{\dal}= (-1)^{\tilde{V}}\bp V^{\dal}$ and
$V^{\dal}$ arbitrary. Moreover, the equivalence class $V\bmod  
\Gamma(M, \Img [\bar{\p}, \ldots])$ always has a representative
of the form 
$$
\sum_{\al=1}^n V^{\al}\frac{\p}{\p z^{\al}}
$$
where $V^{\al}$ are defined uniquely up to addition of a $\bp$-exact term.
The Dolbeault theorem completes the proof. 
\hfill $\Box$

\bip

\bip

\bip

{\bf 2. Deformation theory.} We start with the deformation theory of 
dg-manifolds and then apply the technique to extended complex manifolds.

\bip

{\bf 2.1. Dg-manifolds.} An extended complex manifold is a particular case
of a {\em complex differential $\Z_2$-graded}\, (dg-, for short)  manifold
 which is, by definition, a complex supermanifold equipped with a homological
holomorphic vector field (cf.\ \cite{Ko2}). Morphisms of dg-manifolds, 
$$
(\cX, \p_{\cX}) \stackrel{f}{\lon} (\cS, \p_{\cS}),
$$
are defined as in sect.\ 1.4. It is easy to see that the resulting category
is closed with respect to the fibered products. Note that the fibres of $f$
 are not, in general, dg-manifolds except over the points
where $\p_{\cS}$ vanishes. In this context we define a {\em pointed}\,
dg-manifold as a triple $(\cS, \p_{\cS}, *)$, where $(\cS, \p_{\cS})$
is a (formal) dg-manifold and $*$ a point in $\cS$ such 
that $\p_{\cS}I_{*}\subset  I_{*}^2$, $I_*$ being the ideal sheaf of $*$.

\bip

{\bf 2.1.1. Remark.} According to Kontsevich \cite{Ko2}, any 
Lie superalgebra structure, $ [\ , \ ]$, on a vector superspace
$\fg$ can be equivalently interpreted as a quadratic homological
vector field on $(\Pi\fg, 0)$ viewed as a pointed formal  supermanifold.
Thus $\bH^*(M, T_{M})$ in example~1.5, and $\Pi\bH(\cM,\p_{\cM})$ in
1.4 are  canonically  pointed dg-manifolds.

\bip

{\bf 2.1.2. Lemma.} {\em Let $(\cS, \p_{\cS}, *)$ be a pointed dg-manifold.
Then $\Pi \cT_{\cS, *}$, the tangent space at $*$ with the reversed parity, is 
canonically a Lie superalgebra.}

\sip

\Proof Let $v_1$ and $v_2$ be two tangent vectors at $*\in \cS$, and
$V_1$ and $V_2$ any two germs of holomorphic vector fields 
such that $V_1|_{*}=v_1$ and $V_2|_{*}=v_2$. The bilinear
skew-symmetric operation,
$$
\Ba{rccc}
[\ ,\ ]_*: & \Pi \cT_{ \cS,*}  \times \Pi \cT_{\cS,*} & 
\lon & \Pi T_{\cS,*} \vspace{2mm}
\\
         & \Pi v_1 \times \Pi v_2 &\lon&  \Pi [V_1, [\p_{\cS}, V_2]]\mid_{*}
\Ea
$$
is well defined. The identity $[\p_{\cS},\p_{\cS}]=0$ implies the Jacobi
identity for $[\ , \ ]_*$.
\hfill $\Box$

\bip
{\bf 2.2. Deformations.} A {\em deformation}\, of a complex dg-manifold 
$(\cM,\p_{\cM})$ is, by definition, a morphism, 
$$
(\cX, \p_{\cX}) \stackrel{f}{\lon} (\cS, \p_{\cS},*)
$$
from a complex dg-manifold $(\cX, \p_{\cX})$ to a 
pointed complex dg-base $(\cS, \p_{\cS},*)$ such that
\Bi
\item[(a)] $(f^{-1}(*), \p_{\cX}|_{f^{-1}(*)}) \stackrel{i}{\simeq} 
(\cM, \p_{\cM})$, and \\
\item[(b)] the associated morphism of complex supermanifolds,
$f:\cX \rar \cS$ is locally trivial in an neighbourhood
of $*\in S$.
\Ei

\sip

Given two deformations, 
$f: (\cX, \p_{\cX}) \rar (\cS, \p_{\cS},*)$ and
$\tilde{f}: (\tilde{\cX}, \p_{\tilde{\cX}}) \rar 
(\tilde{\cS}, \p_{\tilde{\cS}},*)$,
of the same dg-manifold manifold $\cM$, a {\em morphism}\, from 
the first to the second, is, by definition, a commutative diagram 
$$
\diagram
&(\cM,\p_{\cM}) \dlto_i \drto^{\tilde{i}} & \\ 
(\cX, \p_{\cX}) \dto_f \rrto^{m}& & (\tilde{\cX}, \p_{\tilde{\cX}}) 
\dto^{\tilde{f}} \\
(\cS, \p_{\cS}, *)\rrto^s && (\tilde{\cS}, \p_{\tilde{\cS}},*),
\enddiagram
$$
where $m$ (resp.\ $s$) is a morphism of (resp.\ pointed) dg-manifolds.
If $s$ is the identity morphism and $m$ is an isomorphism, then the resulting
morphism is called an {\em equivalence}\, of deformations.

\sip

As usually, one defines a deformation of a smooth dg-manifold
over a {\em germ}\, of pointed (possibly, singular) dg-bases. From now on 
we consider only such deformations.

\sip

\bip

{\bf 2.3. Remark.} 
It may look puzzling that we allow in defintion~2.2
the homological vector field $\p_{\cS}$ 
to be non-zero --- the fibres of $f$ over the points where $\p_{\cS}\neq 0$
are {\em not}, in general, dg-manifolds. A more natural definition of 
deformation would be a version of 2.2 with $\p_{\cS}\equiv 0$,
  where the base is simply a pointed analytic
superspace (let us term these {\em $0$-deformations}). 
There are two advantages with understanding a deformation as in 
2.2:
\Bi
\item[(i)] By allowing $\p_{\cS}\neq 0$, we do not loose generality;
$0$-deformation is a deformation. 
\item[(ii)] The $0$-deformation theory of $(\cM,\p_{\cM})$ is, in general,
obstructed, and  its versal moduli
space $(\cV_{\mathsf versal},*)$, when it exists, is a singular analytic space.
The 
role of $\p_{\cS}$ in 2.2 is to overcome 
all these obstructions  producing thereby a {\em smooth}\, versal moduli
dg-space, $(\cS, \p_{\cS}, *)$,  of which 
$\cV_{\mathsf versal}$ is merely a (singular) analytic 
subspace given by the equations $\p_{\cS}=0$ (cf.\ \cite{Me2}). 
In a sense, the deformation theory~2.2 is a smooth resolution
of the more natural $0$-deformation theory.
\Ei

\bip

{\bf 2.4. Proposition.} {\em Let $f: (\cX, \p_{\cX}) \rar (\cS, \p_{\cS},*)$
be a deformation of a dg-manifold $(\cM, \p_{\cM})$. There is a canonical
(even) morphism of Lie superalgebras,}
$$
Df: (\Pi\cT_{\cS,*}\, , [\ , \ ]_*) \lon 
(\bH(\cM, \p_{\cM}), [\ ,\ ]_{\bH}).
$$

\bip

\Proof Let is first construct an  odd linear morphism,
$$
df: \cT_{\cS,*} \lon \bH(\cM, \p_{\cM}),
$$
of vector superspaces. Fixing a local trivialization, $\phi: \cX\simeq 
\cS\times \cM$, we may decompose\footnote{We apologize for being a bit sloppy
in formulating this decomposition.}, 
$$
\p_{\cX}= \p_{\cM} + \p_{\cS} + \Gamma,
$$
for some $f$-vertical odd holomorphic field $\Gamma$ which vanishes
at \mbox{$\cM\times \, *$}. 
Then the homology condition $[\p_{\cX}, \p_{\cX}]=0$ translates into
the Maurer-Cartan(-like) equations,
$$
[\p_{\cM}, \Ga] + [\p_{\cS},\Ga] + \frac{1}{2}[\Ga,\Ga]=0.
$$ 

For a $v$ in $\cT_{\cS,*}$, we define a global holomorphic  vector field
on $\cM$,
$$
df(v):= [V,\Gamma]\mid_{f^{-1}(*)}
$$
where $V$ is an arbitrary extension of $v$ to a germ
of holomorphic vector field on $(\cS,*)$. 
The Maurer-Cartan equations and the fact that $\p_{\cS}$ has  zero
at $*$ of second order imply
$$
[\p_{\cM}, df(v)]=0.
$$
Moreover, the cohomology class of $df(v)$ in $\bH(\cM, \p_{\cM})$
does not depend on the choice of the trivialization $\phi$ so that the 
map $df$ is well-defined.

\sip

Analogously, for any $v_1,v_2\in \cT_{\cS,*}$ we have
$$
\left[V_1, \left[V_2, [\p_{\cM}, \Ga]\right]\right]
+
\left[V_1, \left[V_2, [\p_{\cS}, \Ga]\right]\right]
+
\frac{1}{2}\left[V_1, \left[V_2, [\Ga, \Ga]\right]\right]=0
$$
implying
$$
\left[\left[V_1, [\p_{\cS},V_2]\right], \Gamma\right]\mid_{f^{-1}(*)}
+ \left[[V_1, \Gamma]\mid_{f^{-1}(*)}, [V_2, \Gamma]\mid_{f^{-1}(*)}\right]
=0 \bmod \Img[\p_{\cM}, \ldots ].
$$
Thus $Df:=df\Pi$ is a morphism of Lie superalgebras.
\hfill $\Box$

\bip

{\bf 2.5. Versality.}  If 
$f: (\cX, \p_{\cX}) \rar (\cS, \p_{\cS},*)$ is  a deformation
of a dg-manifold $(\cM, \p_{\cM})$, and
$g:(\tilde{\cS}, \p_{\tilde{\cS}},*)\rar (\cS, \p_{\cS},*)$ 
is a morphism of germs of pointed dg-spaces, then
the fibred product,
$$
(\cX, \p_{\cX})\times_{(\cS, \p_{\cS})}(\tilde{\cS}, \p_{\tilde{\cS}}),
$$  
gives rise to the {\em induced}\, deformation, $g^*(f)$,  of $\cM$ over 
the germ  $(\tilde{\cS}, \p_{\tilde{\cS}},*)$. 

\sip 

A deformation $f: (\cX, \p_{\cX}) \rar (\cS, \p_{\cS},*)$ is called
{\em versal}\, if every other deformation 
$f: (\tilde{\cX}, \p_{\tilde{\cX}}) \rar 
(\tilde{\cS}, \p_{\tilde{\cS}},*)$ of the same extended complex 
manifold is equivalent to the inverse image, $g^*(f)$, of $f$ under some 
morphism of germs, 
$g:(\tilde{\cS}, \p_{\tilde{\cS}},*)\rar (\cS, \p_{\cS},*)$,
 of pointed dg-spaces. A versal deformation $f$ is called {\em minimal}
\, if the associated morphism of Lie superalgebras $Df$ (see 
Proposition~2.4) is an isomorphism. Any two minimal versal
deformations of $(\cM,\p_{\cM})$ are isomorphic.

\bip

{\bf 2.6. Theorem.} {\em Let $(\cM,\p_{\cM})$ be a dg-manifold
with $\dim \bH(\cM,\p_{\cM})<\infty $. Then there exists a smooth minimal
versal deformation of $(\cM,\p_{\cM})$.}

\bip

\Proof This is a special case of the Smoothness Theorem~2.5.6 in \cite{Me2}.
\hfill $\Box$


\bip

{\bf 2.7. Example.} To any (compact) smooth manifold $M$ one may associate
a dg-manifold, 
$$
(\cM= \Pi T_M,\ \p_{\cM}={\mathsf de\ Rham\ differential}).
$$
As $\bH(\cM,\p_{\cM})$ always vanishes (an easy exercise), 
this dg-manifold is rigid, i.e.\ its any deformation
is trivial. This is in accord with the topological nature of the example.

\bip

{\bf 2.7.  Deformation theory of extended complex manifolds.}
This is an obvious modification of the deformation theory of dg-manifolds
with  {\em all}\, the notions and results from 2.1--2.6 holding true.
For example, a {\em deformation}\,  of an extended complex manifold 
$(\cM,\p_{\cM})$ is a proper morphism of complex dg-manifolds,
$$
(\cX, \p_{\cX}) \stackrel{f}{\lon} (\cS, \p_{\cS},*)
$$
such that
\Bi
\item[(a)] for each $t\in {\mathsf Zeros}(\p_{\cS})$ the fibre
$(f^{-1}(t), \p_{\cX}|_{f^{-1}(t)})$ is an extended complex manifold,
\item[(b)] $(f^{-1}(*), \p_{\cX}|_{f^{-1}(*)}) \stackrel{i}{\simeq} 
(\cM, \p_{\cM})$, and \\
\item[(c)] the associated morphism of complex supermanifolds,
$f:\cX \rar \cS$ is locally trivial in an neighbourhood
of $*\in S$.
\Ei

\bip

{\bf 2.7.1. Remarks.} 
(i) In view of 1.1(ii), the 
 local triviality condition 2.7(c) ensures that 
the underlying smooth structure of $\cM$
keeps unchanged upon deformation. It is only the homological vector field that 
undergoes deformation. 

\sip

(ii) The induced morphism, 
$f: {\cH}(\cX, \p_{\cX}) \rar (\cH(\cS, \p_{\cS}),*)$,
may not be locally trivial.

\bip

{\bf 2.7.2. Theorem.} {\em  Let $M$ be a compact complex manifold
and $(\cM, \bar{\p})$ the associated extended complex manifold. Then 
there exists a minimal versal deformation of the latter of the form
$$
\left(\cX, \p_{\cX}\right) 
\stackrel{f}{\lon} \left(\bH^*(M, T_{M}), \p, 0\right).
$$
Moreover,
\Bi
\item[(i)] the homological vector field $\p$ is an invariant of the original
complex manifold $M$. 
\item[(ii)] The restriction of the family $f$ to the analytic
subspace ${\mathsf Zeros}(\p)\cap \bH^1(M, T_{M})$
is equivalent to the classical Kuranishi versal deformation of $M$.

\item[(iii)]
If $M$ is a Calabi-Yau manifold, then 
$\p\equiv 0$ and the base of the above versal deformation can be canonically
identified with the Barannikov-Kontsevich moduli space of (partially) extended
complex structures.
\Ei
}

\bip

\Proof The existence of a minimal versal deformation can be infered directly
from Lemma~1.5.2 and Theorem~2.6. Nevertheless, we will give a  detailed
proof of this statement which makes all other claims, (i)--(iii), almost
obvious. The idea of the proof is to show that, 
though the solution space,
$\widehat{{\mathsf MC}}$,
of Maurer-Cartan equations arising in our extended deformation theory
is much larger than the solution space, ${\mathsf MC}$,  
of the Maurer-Cartan
equations in the original purely algebraic
$(\Gamma(M, T_{M}\ot \Omega^{0,\bullet}_{M}), \bar{\p})$-approach, 
the gauge group turns out to be much larger as well, and, crucially,
at the level of quotients,
$$
\frac{\widehat{\mathsf MC}}{\widehat{\mathsf gauge\ group}}
= \frac{{\mathsf MC}}{{\mathsf gauge\ group}},
$$
we have a canonical isomorphism.

\sip

Let $f: (\cX, \p_{\cX}) \rar (\cS, \p_{\cS},*)$ be a deformation
of $(\cM, \bp)$. 
We may assume without loss of generality that $(\cS, \p_{\cS},*)$ is dual
to a differential Artin superalgebra $(A, \p_A)$, cf.\ \cite{Me2}.
As in the proof of Proposition~2.4, we fix a local trivialization, 
$\phi: \cX\simeq \cS\times \cM$, and decompose$^2$,
$$
\p_{\cX}= \p_{\cM} + \p_{\cS} + \Gamma,
$$
where $\Gamma$ is an $f$-vertical odd holomorphic field $\Gamma$ vanishing
at \mbox{$\cM\times \, *$}. If $\{z^{\al}, z^{\dal}, \psi^{\dal}\}$
is a natural local coordinate system on $\cM$, $\{t^i\}$ a 
local coordinate system on $\cS$ centered at $*$, then $\phi^{-1}$ 
maps the cartesian product
of these into a local coordinate system on $\cX$ in which $\Gamma$ 
can written as follows
$$
\Gamma = \sum_{\al=1}^n \Gamma^{\al}\frac{\p}{\p z^{\al}} + 
\sum_{\dal=1}^n\left(\Gamma^{\dal}
\frac{\p}{\p z^{\dal}} + \gamma^{\dal}\frac{\p}{\p \psi^{\dal}}\right),
$$
for some local functions, $\Gamma^{\al}, \Gamma^{\dal}, \gamma^{\dal}$, of 
$\{z^{\al}, z^{\dal}, \psi^{\dal}, t^i\}$ which vanish at $t^i=0$.
Hence, the  local smooth map,
\begin{eqnarray*}
z^{\al} & \lon & z^{\al},\\
z^{\dal} & \lon & z^{\dal},\\
\psi^{\dal} & \lon & \psi^{\dal} + 
\Gamma^{\dal}(z^{\be}, z^{\dbe}, \psi^{\dbe},t^i),\\
t^i & \lon t^i,
\end{eqnarray*}
is invertable in a small open neighbourhood of $f^{-1}(*)$ in $\cX$.
It is easy to see that this transformation sends the component
$\Gamma^{\dal}\p/\p z^{\dal}$ in $\Gamma$ to zero. Put another way,
this component in $\Gamma$ can always be eliminated by an appropriate
choice --- ``gauge'' ---  of the trivialization $\phi$. 
Then the homology condition
 $[\p_{\cX}, \p_{\cX}]=0$ immediately implies that, in this gauge,
$\gamma^{\dal}=0$ as well. (In physics jargon, both unwanted fields,
$\Gamma^{\dal}$ and $\gamma^{\dal}$,
 correspond to purely gauge degrees of freedom, and, moreover,
 can be eliminated by one single gauge transform as above). 
Of the remaining coordinate transformations preserving the gauge,
only the following ones,
\begin{eqnarray*}
z^{\al} & \lon & z^{\al} + h^{\al}(z^{\be}, z^{\dbe}, \psi^{\dbe},t^i),
\ \ \  h^{\al}\mid_{t=0}=0, \hspace{3cm}  (G)\\
z^{\dal} & \lon & z^{\dal}\\
\psi^{\al} & \lon & \psi^{\dal}\\
t^i & \lon & t^i
\end{eqnarray*}
do effectively change $\Gamma$. Now it is clear
that, in the gauge $\Gamma^{\dal}=\gamma^{\dal}=0$,
 the solution space, $\widehat{\mathsf MC}$, of 
Maurer-Cartan equations in the deformation theory of extended
complex manifolds,
$$
[\bp, \Ga] + [\p_{\cS},\Ga] + \frac{1}{2}[\Ga,\Ga]=0,
$$ 
modulo the remaining gauge freedom $(G)$ can be canonically 
identified with 
the solution space of the Maurer-Cartan equations
of the dg-Lie algebra $(\fg=\Gamma(M, T_{M}\ot \Omega^{0,\bullet}_{M}),
\bar{\p})$,
$$
{\mathsf MC} =
\left\{\Gamma\in (\fg\ot m_A)_{\tln} \mid \bar{\p}\Gamma
+ {\p_A}\Gamma + \frac{1}{2}[\Gamma, \Gamma]=0\right\},
$$
modulo the following gauge transformations,
$$
\Gamma \rar \Gamma^g = e^{\ad_g}\Gamma - 
\frac{e^{ad_g}-1}{\ad_g}(d+{\p}_A)g,
 \ \ \ \ \forall g\in (\fg\ot m_A)_{\tl0}.
$$
Here $m_A$ stands for the maximal ideal in the Artin superalgebra $A$.
It is proved in \cite{Me2} that the associated deformation functor
$$
\Ba{rccc}
{\Df}^{*}_{\fg}: & \left\{\Ba{c} \mbox{the category 
   of differential}\\
                    \mbox{Artin superalgebras}
                    \Ea \right\}&
\lon & \left\{\mbox{the category of sets}\right\}\vspace{3mm}\\
& (A,\p_A) & \lon & {\mathsf MC}/{\mathsf gauge\ group}
\Ea
$$
is always unobstructed, and can be versally represented  by
$\left(\bH^*(M, T_{M}), \p, 0\right)$, where the homological
vector field $\p$ is an invariant of $\fg$. (The 
latter was called in \cite{Me2}
Chen's vector field as its origin  can be traced back
to Chen's power series connection in the theory of iterated integrals.) 

\sip

With the established isomorphism of versal moduli spaces,
the statement 2.7.2(iii) becomes an obvious corollary of
Lemma~2.1 in \cite{BaKo}.

\sip

It remains to check the existence of the classical limit 2.7.2(ii). 
Picking up a Hermitian metric on $M$,
we can construct the adjoint, $\bar{\p}^*$, of the Dolbeault operator
$\bar{\p}$ on $\fg$, the 
Laplacian $\square=\bar{\p}\bar{\p}^*  + \bar{\p}^*\bar{\p}$,
the Green function $G$, and we can identify the cohomology
space $\bH^{\bullet}(M, T_{M})$ with the space of harmonic
elements, $\Ker \bar{\p} \cap \Ker \bar{\p}^*$, in $\fg$. 
Choosing next a harmonic basis, $e_i$, in  
$\bH^{\bullet}(M, T_{M})$, and denoting the associated linear
coordinates by $t^i$ we can represent the Chen's vector
field,  $\p=\sum_i \p^i(t) \p/\p t^i$, as follows \cite{Me2}
$$
\sum_i\p^i(t)e_i = -\frac{1}{2} P[\Gamma, \Gamma],
$$
where $P:\fg \rar \bH^{\bullet}(M, T_{M})$ is the natural projection
to the harmonic constituent, and $\Gamma=\bigoplus_{n=1}^{\infty}\Gamma_{[n]}$
is given by a recursive formula,
\Beqrn
\Gamma_{[1]} &=& \sum_{i} t^{i}e_i 
 \nonumber \\
\Gamma_{[2]} &=& - \frac{1}{2}  G\bar{\p}^* [\Gamma_{[1]}(t),
\Gamma_{[1]}(t)],\nonumber\\
\Gamma_{[3]} &=&  - \frac{1}{2}  G\bar{\p}^*\left([\Gamma_{[1]}(t), \Gamma_{[2]}(t)] +
[\Gamma_{[2]}(t), \Gamma_{[1]}(t)]\right),\nonumber\\
\ldots && \nonumber\\
\Gamma_{[n]} &=& - \frac{1}{2}  G\bar{\p}^* \left(\sum_{k=1}^{n-1}
[\Gamma_{[k]}(t),
\Gamma_{[n-k]}(t)]
\right) \label{split} \\
\ldots && \nonumber
\Eeqrn
It is now obvious that ${\mathsf Zeros}(\p)\cap \bH^1(M,T_{M})$
is precisely the Kuranishi analytic subspace \cite{Ko}. Moreover,
the fibres, $f^{-1}(t)$, of our minimal versal deformation $f$
over this subspace
are precisely the extended complex manifolds, $(\cM_t,\bp_t)$,
 associated, via the 
construction~1.5, to the usual
complex manifolds $M_t$ lying over $t$ in the classical
Kuranishi versal deformation
of $M$.
\hfill $\Box$

\bip

\sip

{\small
{\em Acknowledgement.} This work was done during author's visit
to the Max Planck Institute for Mathematics in Bonn. Excellent working
conditions in the MPIM are gratefully acknowledged.
I would like 
to thank Yu.I.\ Manin for many discussions. }


{\small

{\small
\begin{tabular}{l}
Max Planck Institute for Mathematics in Bonn, and\\
Department of Mathematics, University of Glasgow\\
sm@maths.gla.ac.uk
\end{tabular}
}

\end{document}